\documentclass[11pt,a4paper]{article}
\usepackage[utf8]{inputenc} 
\usepackage{eurosym}

\usepackage{epic,eepic,epsfig,amsmath,amssymb,amsthm}
\usepackage{t1enc}

\usepackage[francais,english]{babel}

\usepackage{caption}

\usepackage{color}

\usepackage{bbm}

\def\virgp{\raise 2pt\hbox{,}}

\renewcommand{\geq}{\geqslant}
\renewcommand{\leq}{\leqslant}
\def\N{{\mathbb N}}

\def\R{{\mathbb R}}

\def\virgp{\raise 2pt\hbox{,}}
\def\cdotpv{\raise 2pt\hbox{;}}

\def\1{\mathbbm{1}}

\newtheorem{theorem}{Theorem}[section]

\newtheorem{proposition}[theorem]{Proposition}
\newtheorem{lemma}[theorem]{Lemma}

\newtheorem{pte}[theorem]{Property}

\theoremstyle{remark}
\newtheorem{remark}{Remark}[section]

\theoremstyle{definition}
\newtheorem{definition}{Definition}[section]

\newtheorem*{notation}{Notation}

\theoremstyle{definition}

\theoremstyle{definition}

\begin{document}

\title{Optimization on fractal sets}

\author{Nizar Riane$^\dag$, Claire David$^\ddag$}

\maketitle
\centerline{$^\dag$ Universit\'e Mohammed V, Agdal, Rabat, Maroc\footnote{nizar.riane@gmail.com} }
\vskip 0.5cm
 
\centerline{$^\ddag$ Sorbonne Universit\'e}

\centerline{CNRS, UMR 7598, Laboratoire Jacques-Louis Lions, 4, place Jussieu 75005, Paris, France\footnote{Claire.David@Sorbonne-Universite.fr}}

\begin{abstract}
We outline necessary and sufficient condition to the existence of extrmas of a function on a self-similar set, and we describe discrete gradient algorithm to find the extrema.
\end{abstract}

	\maketitle
	\vskip 1cm
	

	\vskip 1cm
	
	\vskip 1cm

\section{Introduction}

 \hskip 0.5cm  In 1636, in a correspondence with Martin~Mersenne, Pierre de~Fermat established a necessary condition for the existence of the minimum and the maximum of a function~\cite{Fermat} :

\begin{center}
\textit{``When a quantity, for example the ordinate of a curve, reached its maximum or its minimum, in a situation infinitely close, his increase or decrease is null.''}
\end{center}


 Since then, many results have been set, ranging from free and constraint conditions, giving birth to numerical algorithms that enable one to find the extremas of a function. \\
 
 But, until now, optimization has mainly concerned regular domains, without there being really specific results for fractal sets. The following citation of~G.~Hardy~\cite{Hardy1916}, on a close subject, and even if some may here call it a truism, perfectly reflects the fact that it was ``in consequence of the methods employed''.\\


   Recently, the birth of analysis on fractals, especially, the work of~J.~Kigami~\cite{Kigami1989}, \cite{Kigami1993}, \cite{Kigami1998}, \cite{Kigami2001}, \cite{Kigami2003}, bridged this gap, by giving rise to the building of local operators, which, for a function~$u$ defined on a specific fractal set~$\cal F$, at a given point~\mbox{$X\,\in\, {\cal F}$}, are ``equals to the limit, in a suitable renormalized sense, of the difference between an average value of the function in a neighborhood of~$X$ and~$u(X)$''~\cite{Strichartz1999}. The intrinsinc properties of those operators, analogous to differential ones, make them play the role of Laplacians. If related numerical methods have been developped, for instance, by~R.~Strichartz~\cite{Strichartz1995}, \cite{Strichartz1999}, \cite{Strichartz2000}, \cite{Strichartz2001}, \cite{Strichartz2003}, \cite{Strichartz2012}, \cite{StrichartzLivre2006}, the field of analysis on fractals remains not completely explored, in particular as regards optimization.\\

  In the sequel, in the spirit of Fermat paper, we try to extend results of smooth analysis on extremas in the case of fractals. The novelty of our work lays in the use of the aforementioned differential operators,  specifically designed for fractal sets. To begin with, we examin existence conditions,  then, we present a numerical algorithm that enable us to find local extrema of a continuous function defined on a fractal set.\\

\section{Framework of the sudy}

\hskip 0.5cm In the following, we place ourselves in the euclidian space of dimension~$d\,\in\,\left \lbrace 1,2,3 \right \rbrace $.


\begin{notation}  
	We will classically denote by~$\N^\star$ the set of strictly positive integers.
\end{notation}

\vskip 1cm
\subsection{Self-similar sets}

 \begin{notation}  
In the sequel,~$\cal F$ denotes a fractal domain of Hausdorff dimension~$D_{H}\left({ \cal F}\right)$ ;~$N$ is a strictly positive integer, and~\mbox{$\left \lbrace f_1,\hdots, f_N \right \rbrace $} is a set of contractive maps, where, for any integer~$i$ of~\mbox{$\left \lbrace 1, \hdots, N\right \rbrace $},~\mbox{$R_i\, \in\, ]0,1[ $} is the contraction ratio of~$f_i$, and~\mbox{$P_i\, \in\,\R^d $} the fixed point of~$f_i$.
 	
 \end{notation}
  

 \vskip 1cm

\begin{theorem}\textbf{Gluing Lemma~\cite{Barnsley1985}}\\

	\noindent Given a complete metric space~$(E,\delta)$, a strictly positive integer~$N$, and a set~\mbox{$\left \lbrace f_i\right \rbrace_{1 \leq i \leq N }$}  of contractions on~$E$ with respect to the metric~$\delta$, there exists a unique non-empty compact subset~$K\subset E$ such that:
	
	$$   K =  \underset{  i=1}{\overset{N }{\bigcup}}\,f_{i}\left ( {K} \right )$$
	
	
	\noindent The set~$K $ is said~\textbf{self-similar} with respect to the family~\mbox{$\left \lbrace f_1, \hdots, f_{N } \right \rbrace$}, and called attractor of the iterated function system (IFS)~\mbox{$\left \lbrace f_1 , \hdots, f_{N }  \right \rbrace$}.

\end{theorem} 
 
 \vskip 1cm


\vskip 1cm 

\begin{definition}\textbf{Boundary (or initial) graph}\\
	 
	
	
	
	\noindent We will denote by~$V_0$ the ordered set of the (boundary) points:
	
	$$\left \lbrace P_{1},...,P_{N  }\right \rbrace$$

	\noindent The set of points~$V_0$, where, for any~$i$ of~\mbox{$\left \lbrace  1,...,N -1 \right \rbrace$}, the point~$P_i$ is linked to the point~$P_{i+1}$, constitutes an {oriented graph}, that we will denote by~$  {{\cal F}_0}$.\\
	
	\noindent $V_0$ is called the {set of vertices} of the graph~$ {{\cal F}_0}$. 
	
	\noindent If the point~$P_{N-1}$ is linked to the point~$P_0$, the boundary graph is {closed}.
	
\end{definition}

\vskip 1cm

\begin{definition}\textbf{$m^{th}$ order graph,~$m\,\in\, \N^\star$}\\

\noindent For any strictly positive integer~$m$, we set:
$$V_m =\underset{  i=1}{\overset{N }{\bigcup}}\, f_i \left (V_{m-1}\right ).$$

\noindent The set of points~$V_m$, where two consecutive points are linked, is an oriented graph, which we will denote by~$  {{\cal F}_m}$.\\
\noindent $V_m$ is called the {set of vertices} of the graph~$  {\cal F}_m$. 

\noindent By extension, we will write:

$${\cal F}_m=\underset{  i=1}{\overset{N}{\bigcup}}\, f_i \left ({\cal F}_{m-1}\right ).$$
 
\end{definition}

\vskip 1cm

\begin{pte}
	
	\vskip 0.5cm
	\noindent For any natural integer~$m$:
	$$V_m \subset V_{m+1} .$$
	
\end{pte}

\vskip 1cm

\begin{pte} 
	
	\vskip 0.5cm
	\noindent The set~$\underset{ m\in \N}{ {\bigcup}}\,   V_{m  }  $ is dense in~$  {\cal F}  $.
	
\end{pte}
\vskip 1cm

\begin{definition}\textbf{Word}
	
	\vskip 0.5cm
	
	\noindent Given a strictly positive integer~$m  $, we will call \textbf{number-letter} any integer~${\cal W}_i$ of~\mbox{$\left \lbrace 1, \hdots, N   \right \rbrace $}, and \textbf{word of length~$|{\cal W}|=m$}, on the graph~$ { \cal F}_m $, any set of number-letters of the form:

	$${\cal W}=\left ( {\cal W}_1, \hdots, {\cal W}_m\right).$$
	
	\noindent We will write:
	
	$$f_{\cal W}= f_{{\cal W}_1} \circ \hdots \circ  f_{{\cal W}_m}  .$$

\end{definition}

\vskip 1cm

\begin{definition}\textbf{Vertices }\\
\noindent Two points~$X$ and~$Y$ of~${{ \cal F} }$ will be called \textbf{\emph{vertices}} of the graph~$ { \cal F} $ if there exists a natural integer~$m$ such that:

$$(X,Y) \,\in\, V_m^2$$

\end{definition}

\vskip 1cm
 
\begin{definition}\textbf{Edge relation}\\
	
\noindent Given a natural integer~$m$, two points~$X$ and~$Y$ of~$ {{ \cal F}_m}$ will be called \emph{\textbf{adjacent}} if and only if~$X$ and~$Y$ are neighbors in~$ {{ \cal F}_m}$. We will write:

$$X \underset{m }{\sim}  Y$$

\noindent This edge relation ensures the existence of a word~{${\cal W}=\left ( {\cal W}_1, \hdots, {\cal W}_m\right)$} of length~$ m$, such that~$X$ and~$Y$ both belong to the iterate:

$$f_{\cal W} \,V_0=\left (f_{{\cal W}_1} \circ \hdots \circ  f_{{\cal W}_m} \right) \,V_0$$

\noindent Given two points~$X$ and~$Y$ of~$ { \cal F} $, we will say that~$X$ and~$Y$ are \textbf{\emph{adjacent}} if and only if there exists a natural integer~$m$ such that:
$$X  \underset{m }{\sim}  Y$$
 \end{definition}

\vskip 1cm

\begin{definition}\textbf{Adresses}\\
 
	Given a natural integer~$m$, and a vertex~$X$ of~${\cal F}_m$, we will call {address of the vertex~$X$} an expression of the form

	$$X= f_{ \cal W} \left (P_i\right ) $$

\noindent	where~$\cal W$ is a word of length~$m$, and~$i$ a natural integer in~\mbox{$\left \lbrace 1,\hdots, N \right \rbrace $}.
\end{definition}

\vskip 1cm

\begin{pte}\textbf{Subcell - Junction points}\\

\noindent Given a natural integer~$m$, the graph~${\cal F}_m$ can be written as the finite union of~$N^m$ subgraphs:

$${\cal F}_m =  \underset{  |  {\cal W} |  =m}{\overset{ }{\bigcup}}\, f_{ \cal W} \left ({ \cal F}_0\right )$$

\noindent For any word~$\cal W$ of length~$m$,~$f_{ \cal W} \left ({ \cal F}_0\right )$ will be called~{$m^{th}$-order cell}, or {subcell}.

\end{pte}

\vskip 1cm

\begin{notation}

\begin{itemize}

	\item[\emph{i}.] Given a strictly positive integer~$m$, we will denote by $\Sigma_m$ to be the set of words~\mbox{${\cal W} \, \in \,  \{1,\hdots ,N\}^m$} of length $m$.\\

		\item[\emph{ii}.] We then set:
			
			$$\Sigma=\underset{ m \in\N^\star }{{\bigcup}}\, \Sigma_m$$

\end{itemize}
\end{notation}

\vskip 1cm

\begin{notation}For the sake of clarity, we will, from now on, refer to a self-similar set either by~$\cal F$ or by:

	$$\left({\cal F},\Sigma,\left(f_i\right)_{i=1}^N\right)$$

\end{notation}

\vskip 1cm

\begin{notation}
	We will denote by:

\begin{itemize}

	\item[\emph{i}.]  $\sigma$, the~\textbf{shift map} from~$\Sigma$ to~$\Sigma$ which, for any word $\cal W$, deletes the first \textbf{"letter"} i.e.:

$$ \sigma(12233\hdots )=2233 \hdots$$

\item[\emph{ii}.]  $\pi$, the surjective map from~$\Sigma$ to $\cal F$ defined for every infinite "word" \mbox{${\cal W}={ \cal W}_1 { \cal W}_2 \hdots \, \in\, \Sigma$}, by:

$$ \pi (W) =\underset{ m \in\N^\star }{{\bigcap}}\,   f_{{ \cal W}_1 \hdots  { \cal W}_m}({ \cal F}) $$

\noindent where $f_{{ \cal W}_1 \hdots { \cal W}_m}=f_{{ \cal W}_1}\circ  \hdots\circ f_{{ \cal W}_m}$.\\


\item[\emph{iii}.] $$ C_{L,{ \cal F}}= \underset{(i,j)\,\in\, \Sigma^2, \, i \neq j}{{\bigcup}}\,   \left  ({ \cal F}_i \cap { \cal F}_j \right )$$

\end{itemize}
\end{notation}

\vskip 1cm

\begin{definition}

\noindent We define:

\begin{itemize}

	\item[\emph{i}.] The critical set:

$$ C_L = \pi^{-1}\left (C_{L,{\cal F}}\right )$$

\item[\emph{ii}.] The post-critical set:

$$ P = \underset{ m \in\N^\star }{{\bigcup}}\,  \sigma^{m}(C_{L})$$

\end{itemize}
	
	\end{definition}
\vskip 1cm

\begin{pte}
$$V_0=\pi \left (  \underset{ m \in\N^\star }{{\bigcup}}\,  \sigma^{m}(C_{L})  \right )$$


\end{pte}
  
\vskip 1cm

\begin{notation}
	\noindent Given a natural integer~$m$, we will denote by~$\mathcal{N}_m$ the number of vertices of the graph ${\cal F}_m$, {$\mathcal{A}_m$ the edge set of ${\cal F}_m$}, and by~$C$ the cardinal of the set~$  C_{L,{\cal F}}$.
	
	
	\end{notation}

		\vskip 1cm

\begin{proposition}
\noindent  One has:

$$\mathcal{N}_0  =N_0$$
		
\noindent and, for any strictly positive integer~$m$:

$$
\mathcal{N}_m    =N\times \mathcal{N}_{m-1}- C \quad , \quad 
\# A_m  = N\times \# A_{m-1}$$ 
\end{proposition}
\vskip 1cm

\begin{proof}
The graph ${\cal F}_m$ is the union of $N$ copies of the graph ${\cal F}_m$. Each copy has~$\#A_{m-1}$ edges, and shares vertices with others copies. One can thus consider the copies as the vertices of a complete graph which has a number of edges equal to $C$, so there are~$C$ vertices to discount.
\end{proof}

\vskip 1cm

\begin{remark}
\noindent One may check that:

$$A_m=N^m A_0=\displaystyle{N^m \, (N_0(N_0-1))} \quad , \quad \mathcal{N}_m=\displaystyle{N^m N_0-\left(\frac{1-N^m}{1-N}\right)C}=\mathcal{O}(N^m).$$
\end{remark}

\vskip 1cm

\begin{definition}{\textbf{System of neighborhood~\cite{Kigami2001} \\}}
	
\noindent Let~$\left({\cal F},\Sigma,\left(f_i\right)_{i=1}^N\right)$ be a self-similar structure. For any~\mbox{$X\,\in \,{\cal F}$}, and any natural integer~$m$, we set:

$$ {\cal F}_{m,X} = \underset{ {\cal W}\in \Sigma_m, \, X\, \in\,  f_{\cal W}({\cal F}) }{{\bigcup}}\, 
 f_{\cal W}({\cal F})$$

\noindent which will be called system of neighborhood of $X$.
\end{definition}
\vskip 1cm


\vskip 1cm


\begin{definition}{\textbf{Self-similar measure on~$\cal F$~\cite{StrichartzLivre2006}}}\\

	\noindent A measure~$\mu$ with full support on~$\R^d$ will be called \textbf{self-similar measure} on~\mbox{${\cal F}=
		\underset{  i=1}{\overset{N }{\bigcup}}\,f_{i}\left ( {\cal F} \right )$} if, given  a family of strictly positive weights~\mbox{$\left (\mu_i\right)_{1 \leq i \leq N }$} such that:
	
	$$ \displaystyle \sum_{i=1}^{N  } \mu_i=1$$
	\noindent one has:
	$$ \mu= \displaystyle \sum_{i=1}^{N  } \mu_i\,\mu\circ f_i^{-1} $$

\end{definition}

\vskip 1cm
	
\begin{pte}\textbf{Building of a self-similar measure on~$\cal F$}\\
		
\noindent We set, for any integer~$i$ belonging to~$\left \lbrace 1, \hdots, N \right \rbrace$:
		
$$\mu_i=R_i^{D_{H}\left( {\cal F }\right)}$$
		
\noindent  One has then:

$$\displaystyle \sum_{i=1}^N R_i^{D_{H}\left({ \cal F} \right)}=1.$$

\noindent which enables us to define a self-similar measure~$\mu$ on~$\cal F$ through:
		
$$\mu =\displaystyle \,\sum_{i=1}^N \mu_i  \, \mu \circ f_i^{-1}  $$
		
\noindent One may note that the measure~$\mu$ corresponds to the normalized  $D_{H}\left({ \cal F }\right)$-dimensional Hausdorff measure ({$\mathcal{H}^{D_{H}\left({ \cal F }\right)}$})	(we refer to~\cite{Falconer1985}):

$$ \mu(E)=\frac{\mathcal{H}^{D_{H}\left({ \cal F }\right)}(E\cap { \cal F })}{\mathcal{H}^{D_{H}\left({ \cal F }\right)}({ \cal F })}$$

\noindent for any subset $E\subset \R^d$.

\end{pte}

	\vskip 1cm

\subsection{Laplacians, on self-similar sets}


\begin{definition}{\textbf{Energy, on the graph~${\cal F }_m$,~$m \,\in\,\N$, of a pair of functions}}\\

	\noindent Given a natural integer~$m$, and two real valued functions~$u$ and~$v$, defined on the set~$V_m$ of the~vertices of~${\cal F }_m$, we introduce \textbf{the energy, on the graph~${\cal F }_m$, of the pair of functions~$(u,v)$}, as:
	
	$$\begin{array}{ccc}
	{\cal{E}}_{{\cal F }_m }(u,v)
	&= & \displaystyle \sum_{(X,Y)\,\in\,V_m^2,\, X \sim Y}  \left ( u \left (X \right)-u \left (Y\right) \right )\,
	\left ( v \left (X\right)-v \left (Y \right)\right )
	\end{array}
	$$
	
	\noindent For the sake of simplicity, we will write it under the form:
	
	$$ {\cal{E}}_{{\cal F }_m}(u,v)= \displaystyle \sum_{X  \underset{m }{\sim}  Y}\left (u(X)-u(Y)\right )\,\left(v(X)-v(Y)\right) $$

\end{definition}

\vskip 1cm

\begin{definition}\textbf{Dirichlet form on a measured space}\\
 
\noindent (we refer to the paper \cite{Beurling1985}, or the book \cite{Fukushima1994})\\

\noindent Given a measured space~$(E, \mu)$, a \emph{\textbf{Dirichlet form}} on~$E$ is a bilinear symmetric form, that we will denote by~${\cal E}$,
defined on a vectorial subspace~$D$ dense in $L^2(E, \mu) $, such that:\\

\begin{enumerate}
	
	\item For any real-valued function~$u$ defined on~$D$ :  ${\cal E}(u,u) \geq 0$.
	
	\item   $D$, equipped with the inner product which, to any pair~$(u,v)$ of~$D \times D $, associates:

	$$  (u,v)_{\cal{E}}  = (u,v)_{L^2(E,\mu)} + {\cal {E}}(u,v)$$
	
	is a Hilbert space.

	\item For any real-valued function~$u$ defined on~$D$, if:
	$$ u_\star = \min\, (\max(u, 0) , 1) \,\in \,D$$

	\noindent then : ${ \cal{E}}(u_\star,u_\star)\leq { \cal{E}}(u,u)$ (Markov property, or lack of memory property).

\end{enumerate}

\end{definition}

\vskip 1cm

\begin{definition}{Dirichlet form, on a finite set (see~\cite{Kigami2003})}\\
 
\noindent Let~$V$ denote a finite set~$V$, equipped with the usual inner product which, to any pair~$(u,v)$ of functions defined on~$V$, associates:

$$(u,v)= \displaystyle \sum_{p\in  V} u(p)\,v(p)$$

\noindent A \emph{\textbf{Dirichlet form}} on~$V$ is a symmetric bilinear form~${\cal E}$, such that:\\

\begin{enumerate}

\item For any real valued function~$u$ defined on~$V$:  ${\cal E}(u,u) \geq 0$.

\item   $  {\cal {E}}(u,u)= 0$ if and only if~$u$ is constant on~$V$.

\item For any real-valued function~$u$ defined on~$V$, if:~$ u_\star = \min\, (\max(u, 0) , 1)  $, i.e. :

$$\forall \,p \,\in\,V \, : \quad u_\star(p)= \left \lbrace \begin{array}{ccc} 1 & \text{if}& u(p) \geq 1 \\u(p) & \text{si}& 0 <u(p) < 1 \\0  & \text{if}& u(p) \leq 0 \end{array} \right.$$

\noindent then: ${ \cal{E}}(u_\star,u_\star)\leq { \cal{E}}(u,u)$ (Markov property).

\end{enumerate}

\end{definition}
 
\vskip 1cm

\begin{pte}

	\vskip 0.5cm
	\noindent Given a natural integer~$m$, and a real-valued function~$u$, defined on the set~$V_m$ of vertices of~$  {\cal F}_m$, the map, which, to any pair of real-valued, continuous functions~$(u,v)$ defined on~$V_m $, associates:
	$$ {\cal{E}}_{ {\cal F}_m}(u,v)= \displaystyle \sum_{X  \underset{m }{\sim}  Y}\left (u(X)-u(Y)\right )\, \left (v(X)-v(Y)\right ) $$

	\noindent is a Dirichlet form on~$ {\cal F }_m $.\\
	\noindent Moreover:
	
	$${\cal{E}}_{{\cal F }_m }(u,u)=0 \Leftrightarrow u  \text{ is constant}$$

\end{pte}

\vskip 1cm
\begin{proposition}  
	
	\noindent For any strictly positive integer~$m$, if~$u$ is a real-valued function defined on~$V_{m-1}$, its \textbf{harmonic extension}, denoted by~$ \tilde{u}$, is obtained as the extension of~$u$ to~$V_m$ which minimizes the energy:

	$$  {\cal{E}}_{{\cal F}_m}(\tilde{u},\tilde{u})=\displaystyle \sum_{X \underset{m }{\sim} Y} \displaystyle  (\tilde{u}(X)-\tilde{u}(Y))^2 $$

\end{proposition}

\vskip 1cm

\begin{remark} \textbf{Concretely:}\\

\noindent The link between~$   {\cal{E}}_{{\cal F }_m}$ and~$  {\cal{E}}_{ {\cal F}_{m-1}}$ is obtained through the introduction of two strictly positive constants~$r_m$ and~$r_{m-1}$ such that:

$$  r_{m }\, \displaystyle \sum_{X \underset{m  }{\sim} Y} \displaystyle  (\tilde{u}(X)-\tilde{u}(Y))^2  =   r_{m-1} \,\displaystyle  \sum_{X \underset{m-1 }{\sim} Y} (u(X)-u(Y))^2 $$

\noindent In particular:

$$  r_{1 }\, \displaystyle \sum_{X \underset{1  }{\sim} Y} (\tilde{u}(X)-\tilde{u}(Y))^2 =  r_{0}\,\displaystyle  \sum_{X \underset{0 }{\sim} Y} (u(X)-u(Y))^2$$

\noindent For the sake of simplicity, one fixes the value of the initial constant:~$r_0=1$. Then:

$$ {\cal{E}}_{ {\cal F}_1 }(\tilde{u},\tilde{u})= \displaystyle \frac{1}{ r_{1 } }\,  {\cal{E}}_{{\cal F}_0}(\tilde{u},\tilde{u})$$

\noindent We set:

$$r = \displaystyle \frac{1}{r_{1 }} $$

\noindent and:

$$  {\cal{E}}_{m}(u)=r_m\, \sum_{X \underset{m }{\sim} Y}  (\tilde{u}(X)-\tilde{u}(Y))^2 $$

\noindent Since the determination of the harmonic extension of a function appears to be a local problem, on the graph~${\cal F }_{m-1}$, which is linked to the graph~${\cal F }_m$ by a similar process as the one that links~${\cal F }_1$ to~${\cal F }_0$, one deduces, for any strictly positive integer~$m$:

$$ {\cal{E}}_{ {\cal F}_m }(\tilde{u},\tilde{u})= \displaystyle \frac{1}{ r_{1 }}\,  {\cal{E}}_{ {\cal F}_{m-1}}(\tilde{u},\tilde{u})$$

\noindent By induction, one gets:

$$r_m=r_1^m =r^{-m}  $$

\noindent If~$v$ is a real-valued function, defined on~$V_{m-1}$, of harmonic extension~$ \tilde{v}$, we will write:

$$  {\cal{E}}_{m}(u,v)=r^{-m}\, \sum_{X \underset{m }{\sim} Y} 
(\tilde{u}(X)-\tilde{u}(Y)) \, (\tilde{v}(X)-\tilde{v}(Y)) $$

\noindent For further precision on the construction and existence of harmonic extensions, we refer to~\emph{\cite{Sabot1987}}.

\end{remark}

\vskip 1cm

\begin{definition}\textbf{Renormalized energy, for a continuous function~$u$, defined on~$  {\cal F }_m   $,~$m\,\in\,\N$}\\
	 
	\noindent Given a natural integer~$m$, one defines the \textbf{normalized energy}, for a continuous function~$u$, defined on~$  {\cal F }_m   $, by:

	$$
	{\cal{E}}_m (u)=    \displaystyle \sum_{X  \underset{m }{\sim}  Y} r^{-m}\, 
	\left (u (X)-u(Y)\right )^2 $$

\end{definition}

\vskip 1cm

\begin{definition}\textbf{Normalized energy, for a continuous function~$u$, defined on~$  {\cal F }   $} \\
	
	Given a function~$u$ defined on
	$$\displaystyle {V_\star =\underset{{i\in \N}}\bigcup \,V_i}$$

	\noindent one defines the \textbf{normalized energy}:

	$$
	{\cal{E}} (u)=   
	\displaystyle \lim_{m \to + \infty}\displaystyle \sum_{X  \underset{m }{\sim}  Y} r^{-m}\, 
	\left (u (X)-u (Y)\right )^2 $$

\end{definition}
\vskip 1cm

\begin{definition}\textbf{Dirichlet form, for a pair of continuous functions defined on~$ {\cal F }  $} \\

	\noindent We define the Dirichlet form~$\cal{E}$ which, to any pair of real-valued, continuous functions~$(u,v)$ defined on the graph~$  {\cal F }   $, associates, subject to its existence:
	
	$$
	{\cal{E}} (u,v)=   
	\displaystyle \lim_{m \to + \infty}\displaystyle \sum_{X  \underset{m }{\sim}  Y} r^{-m}\, 
	\left (u_{\mid V_m}(X)-u_{\mid V_m}(Y)\right )\,\left(v_{\mid V_m}(X)-v_{\mid V_m}(Y)\right)
	$$

\end{definition}

\vskip 1cm

\begin{notation}
 
\noindent  We will denote by:

\begin{enumerate}
\item[\emph{i}.]  ~$\text{dom}\,{\cal E}$ the subspace of continuous functions defined on~$ {\cal F}$, such that:

$$\mathcal{E}(u)< + \infty$$

 \item[\emph{ii}.] ~$\text{dom}_0\,{ \cal E}$ the  subspace of continuous functions defined on~$ {\cal F}$, which take the value zero on~$V_0$, and such that:
	
	$$\mathcal{E}(u)< + \infty$$
	
\end{enumerate}
	
\end{notation}
\vskip 1cm

\begin{lemma} 

The map: 
$$\begin{array}{ccc}
{\text{dom} \,{ \cal E}}  \big / {\text{Constants}}  \times {\text{dom} \,{ \cal E}}  \big / {\text{Constants}} &  \to& \R  \\
(u,v) & \mapsto & 
{\cal{E}} (u,v) \end{array}$$

\normalsize	
\noindent defines an inner product on~\mbox{$	{\text{dom} \,{ \cal E}}  \big / {\text{Constants}} $}.

\end{lemma}

\vskip 1cm

\begin{theorem} 
	$\left ( {\text{dom} \,{ \cal E}}  \big / {\text{Constants}} , { \cal E} ( \cdot, \cdot ) \right) $ is a complete Hilbert space.
	
\end{theorem}

\vskip 1cm
 
 \begin{definition}\textbf{Harmonic function}\\
 	 
 	\noindent A real-valued function~$u$, defined on~$\displaystyle {V_\star =\underset{{i\in \N}}\bigcup \,V_i}$, will be said to be \textbf{harmonic} if, for any natural integer~$m$, its restriction~$u_{\vert V_m}$ is harmonic:

 	$$\forall\, m\,\in\,\N,\,\forall\, X\,\in\, V_m\setminus V_0 \, :\quad \Delta_m u_{\vert V_m}(X) =0 $$
 	
 \end{definition}
 
 \vskip 1cm

\begin{notation} 

We will denote by~$\text{dom}\,{\cal E}$ the subspace of continuous functions~$u$ defined on~$ {\cal F}$, such that:

$$\mathcal{E}(u)=
\displaystyle \lim_{m \to + \infty}\displaystyle \sum_{X  \underset{m }{\sim}  Y} r^{-m}\, 
\left (u_{\mid V_m}(X)-u_{\mid V_m}(Y)\right )^2 < + \infty$$

\end{notation}

\vskip 1cm
\begin{definition} 
 
\label{Lapl} 
\noindent We will denote by~$\text{dom}\, \Delta$ the existence domain of the Laplacian, on~$ {\cal F}$, as the set of functions~$u$ of~$\text{dom}\, \mathcal{E}$ such that there exists a continuous function on~$ {\cal F}$, denoted by~$\Delta \,u$, that we will call \textbf{Laplacian of~$u$}, such that,~\mbox{$\text{for any } v \,\in \,\text{dom}\, \mathcal{E}\, , \, v_{\vert {\cal F}_0}=0 \, $}:

$$\begin{array}{ccc}\mathcal{E}(u,v)& =&
\displaystyle \lim_{m \to + \infty}\displaystyle \sum_{X  \underset{m }{\sim}  Y} r^{-m}\, 
\left (u_{\mid V_m}(X)-u_{\mid V_m}(Y)\right )\,\left(v_{\mid V_m}(X)-v_{\mid V_m}(Y)\right)\\
&=&-\displaystyle \int_{\cal F}  v\, \Delta u   \,d\mu \end{array}$$

\end{definition}

\vskip 1cm

\begin{theorem}

\vskip 0.5cm

$$\left (u \, \in\,  \text{dom}\, \Delta \quad \text{and} \quad \Delta \, u=0\right) 
\quad \text{\underline{if and only if}~$u$ is harmonic}$$

\end{theorem}

\vskip 1cm

\begin{notation}
	\vskip 0.5cm
	
	\emph{i}. Given a natural integer~$m$,~${\cal S} \left ({\cal H}_0,V_m \right)$ denotes the space of spline functions "of level~$m$",~$u$, defined on~$ {\cal F}$, continuous, such that, for any word~$\cal M$ of length~$m$,~\mbox{$u \circ T_{\cal M}$} is harmonic, i.e.:
	
	$$\Delta_m \, \left ( u \circ T_{\cal M} \right)=0$$

	\emph{ii}.~${\cal H}_0\subset \text{dom}\, \Delta$ denotes the space of harmonic functions, i.e. the space of functions~$u \,\in\,\ \text{dom}\, \Delta$ such that:
	
	$$\Delta\,u=0$$
	
\end{notation}

\vskip 1cm

\begin{pte}
	 
	For any natural integer~$m$:~${\cal S} \left ({\cal H}_0,V_m \right )\subset  \text{dom }{ \cal E}$.
	
\end{pte}

\vskip 1cm

\begin{theorem}\textbf{Pointwise formula}\\
 
Let~$m$ be a strictly positive integer,~$X \,\in \,V_\star \setminus V_0$, and~\mbox{$\psi_X^{m}\,\in\,{\cal S} \left ({\cal H}_0,V_m \right)$} a spline  function such that:

$$\psi_X^{m}(Y)=\left \lbrace \begin{array}{ccc}\delta_{XY} & \forall& Y\,\in \,V_m \\
0 & \forall& Y\,\notin \,V_m \end{array} \right. \quad,  \quad \text{where} \quad    \delta_{XY} =\left \lbrace \begin{array}{ccc}1& \text{if} & X=Y\\ 0& \text{else} &  \end{array} \right.$$

 \begin{enumerate}
\item[ {i}.] For any function~$u$ of~$\text{dom}\, \Delta $, such that its Laplacian exists, the sequence

$$\left ( r^{-m}\, \left \lbrace \int_{\cal F} \psi_X^{m}\,d \mu\right \rbrace^{-1} \,\Delta_m u(X) \right)_{m\in\N}$$

converges uniformly towards

$$  \Delta  u(X) $$

\item[ {ii}.]  Conversely, given a continuous function~$u$ on~$\cal F$ such that the sequence

$$\left ( r^{-m}\, \left \lbrace \int_{\cal F} \psi_X^{m}\,d \mu\right \rbrace^{-1} \,\Delta_m u(X) \right)_{m\in\N}$$

converges uniformly towards a continuous function on~$V_\star \setminus V_0$, one has:

$$  u\, \in\, \text{dom}\, \Delta \quad \text{and} \quad 
\Delta \, u(X)=\displaystyle \lim_{m \to + \infty}  r^{-m}\, \left \lbrace \int_{\cal F} \psi_X^{m}\,d \mu\right \rbrace^{-1} \,\Delta_m u(X)  $$

\end{enumerate}

\end{theorem}

\vskip 1cm

	
		
	
	
	
		
	
	
\subsection{Existence of extrema}

\vskip 1cm

\begin{definition}{\textbf{Extrema }}\\
	
\noindent Given a continuous function~$u$ defined on the fractal set~$\cal F$, and~$X\, \in \, {\cal F}$, we will say that~$u$:
\begin{enumerate}
\item[\emph{i}.] has a global minimum (resp. a global maximum) at~$X$ if:
$$ \forall \, Y\, \in\,  {\cal F} \, :\quad u(X)\leq u(Y) \quad (\text{resp. } u(X)\geq u(Y) ) $$
\item[\emph{ii}.]  a local minimum (resp. a local maximum) at~$X$ if there exist a neighborhood~$V$ of $X$ such that
$$ \forall \, Y\, \in \, V \,:\quad u(X)\leq u(Y) \quad (\text{resp. } u(X)\geq u(Y) ).$$
\end{enumerate}
\end{definition}

\vskip 1cm

\begin{theorem}
	Given a continuous function~$u$ defined on the compact fractal set
	
$${\cal F} =	 \underset{  i=1}{\overset{N }{\bigcup}}\,f_{i}\left ( {{\cal F}} \right )$$

\noindent the Weierstrass extreme value theorem ensures the existence of:

$$\displaystyle  \min_{X\in {\cal F}} \, u  \quad \text{and} \quad \max_{X\in {\cal F}} \, u$$

\end{theorem}

\vskip 1cm

	








\vskip 1cm

\begin{theorem}\textbf{Laplacian test for fractals~\cite{Kigami2001}}\\
	
\noindent Given a continuous, real-valued function~$u$ defined on~$\cal F$, and belonging to $dom_{\Delta }$:

\begin{enumerate}
\item[\emph{i}.] If~$u$ admits a local maximum at~\mbox{$X_0\, \in\, {\cal F}$}, then:
$$ \Delta u(X_0) \leq 0$$
\item[\emph{ii}.] IIf~$u$ admit a local minimum at~\mbox{$X_0\, \in\, {\cal F}$}, then:
$$\Delta u(X_0) \geq 0$$
\end{enumerate}
\end{theorem}

\vskip 1cm

\begin{proof}$\,$\\
\noindent \emph{i}. If $u$ admits a local maximum at $X_0$, then for sufficiently large values of the integer~$m$:
$$ \Delta_{m} u(X_0) = \displaystyle \sum_{X_0  \underset{m }{\sim}  Y} \left(u(Y) - u(X_0)\right) \leq 0.$$

\noindent This enables us to conclude that:

$$ \Delta_{\mu} u(X)=\displaystyle \lim_{m\rightarrow\infty} r^{-m}\,  \left(\int_K \psi^{(m)}_{X_m} d\mu\right)^{-1}\, \Delta_m u(X_m) \leq 0$$

\noindent \emph{ii}. can be proved in a similar way.
 
\end{proof}

\vskip 1cm

\section{Numerical algorithm and dynamic programming}

 \hskip 0.5cm  In the following, we will present a numerical algorithm, based on discrete gradient, to find a local maximizer (resp. local minimizer) of a function continuous $u$ on $K$.\\

\subsection{The algorithm}

 \hskip 0.5cm  We recall that ${\cal F}_m=(V_m,A_m)$ is the oriented graph approximation of $\cal F$ of order $m$, where $V_m$ is the vertices set and $A_m$ is the edge set. We can check that the distance between two connected vertices is of order $2^{-m}$.\\

  In order to find the local maximumizer, we will provide every edge~$\{XY\}$ with the weight $D_{XY}=u(Y)-u(X)$. In order to find an appropriate approximation of the maximizer $X^{\star}$, we fix a degree of tolerance $\varepsilon >0$, and the graph ${\cal F}_m$ of order $m$ such that

$$ 2^{-m} \leq \varepsilon $$

  Starting at an arbitrary point $X_0$ in $V_m$, we follow the direction of the maximal positive gradient at $X_0$, i.e.

$$ \max_{X_0\underset{m}{\sim}Y} \{D_{X_0 Y} \,| \, D_{X_0 Y} >0\}$$

  We replace the initial point by $\arg\max_{X_0\underset{m}{\sim}Y} \{D_{X_0 Y} \,| \, D_{X_0 Y} >0\}$, and we do the same operation until

$$ \max_{X_0\underset{m}{\sim}Y} \{D_{X_0 Y} \,| \, D_{X_0 Y} \geq 0\} = 0$$

  In this case, the algorithm stop and we have the approximation of $X^{\star}$. We can can resume the algorithm in the following steps :\\

\vskip 1cm

\noindent \textbf{Discrete gradient algorithm}

\begin{enumerate}
\item Fix a degree of tolerance~$\varepsilon>0$.
\item Build the graph ${\cal F}_m$, for $m$ such that:
$$ 2^{-m} \leq \varepsilon $$
\item Fix~$X=X_0$ for $X_0\,  \in \, V_m$.
\item While~$\displaystyle \max_{X\underset{m}{\sim}Y} \, \{D_{X Y} \,| \, D_{XY} \geq 0\} > 0$:\\
\centerline{Update~$X=\displaystyle\arg\max_{X\underset{m}{\sim}Y} \{D_{X Y} \,| \, D_{X Y} >0\}$.}
 
\item Return~$X$.
\end{enumerate}

\vskip 1cm

\subsection{Numerical analysis and dynamic programming}

 \hskip 0.5cm The algorithm presented above can be viewed as a dynamical programming algorithm on the directed graph~${\cal F}_m$,~$m\,\in\, \N^\star$~(we refer to~\cite{Gaubert2015}). Values of the function can be calculated recursively:

$$
v^0_X =0 \quad , \quad \forall \, X\, \in \, V_m \quad
v^m_X  =\mathcal{B} \left( v^{m-1}\right)_X =\sup_{Y\in V_m} \left( D_{XY} + v^{n-1}_X \right) \ $$
\noindent where $\mathcal{B}$ denotes the Bellman operator:

\begin{align*}
\mathcal{B} &: \bar{\R}^{V_m} \rightarrow \bar{\R}^{V_m}\\
\mathcal{B}\left( v\right)_X &=\sup_{Y\in V_m} \left( D_{XY} + v_X \right)
\end{align*}

\noindent where $D_{XY}$ is the weight associated to the arrow ${XY}\in A_m$ and $D_{XY} = -\infty$ if ${XY}\not\in A_m $.\\

  As shown in the first section, $\# A_m = \mathcal{O}(N^{m})$ and $\# V_m=\mathcal{O}(N^m)$. The number of possible transitions is at most of order $\# A_m = \mathcal{O}(N^{m})$. Since 
$$\# V_m=\mathcal{O}(N^m)$$

\noindent one can deduce that the calculation time of the maximum is of order~$\mathcal{O}(N^{2\, m})$.\\

\subsubsection{Sierpi\'{n}ski simplices}

 \hskip 0.5cm Sierpi\'{n}ski simplices are sparse graph. Thus, the calculation time required for the gradient algorithm can be optimized: each vertex has a finite number of neighbors, which ensures a calculation time at step~\mbox{$m\,\in\,\N^\star$} which is of order $\mathcal{O}(N^{m})$.\\

  Moreover, computations can be simplified using the fact that every vertex~$X$ has two adresses:

$$ X=f_{{\cal W}_1}(P_i)=f_{{\cal W}_2}(P_j)$$

\noindent where~\mbox{$\left ( {\cal W}_1,  {\cal W}_2\right) \, \in\,  \Sigma^2$} and~\mbox{$\left (P_i,\,P_j\right) \, \in \, V_0^2$}, with~$i\neq j$. Thus, the neighbors of~$X$ are given by:

$$ \left \lbrace  \underset{ k \neq i}{\overset{ }{\bigcup}}\,f_{{\cal W}_2}(P_k) \right \rbrace \, \bigcup \, \left  \lbrace 
 \underset{ \ell \neq j}{\overset{ }{\bigcup}}\,  f_{{\cal W}_2}(P_\ell) \right \rbrace$$

\vskip 1cm
\paragraph{The Sierpi\'{n}ski Gasket\\ \\ } 

\vskip 1cm

  In the case of Sierpi\'{n}ski Gasket, one can optimize the calculation time of the gradient algorithm.  Since, for any natural integer~$m$:

$$\#V_m=\displaystyle\frac{3^{m+1}+3}{2}\quad \text{and} \quad \# A_m=2\times 3^{m+1}$$

\noindent and given the fact that every vertex~\mbox{$X\, \in \, V_m\setminus V_0$} has four neighbors, it follows that the calculation time of the maximum at step~\mbox{$m\,\in\,\N^\star$}  is thus of order~$\mathcal{O}(3^{m})$.\\

  In the sequel (see~figures~1 to~4), we present results of our algorithm in the case of Sierpi\'{n}ski Gasket with vertices:
$$P_0=(0,0) \quad , \quad P_1=(1,0)\quad , \quad P_2=\displaystyle \left(\frac{1}{2},\frac{\sqrt{3}}{2}\right)$$

\noindent  for the value~$m=6$.\\

  The color function is related to the gradient of the one at stake, high values ranging from red to blue.  

\begin{center}
\includegraphics[scale=1]{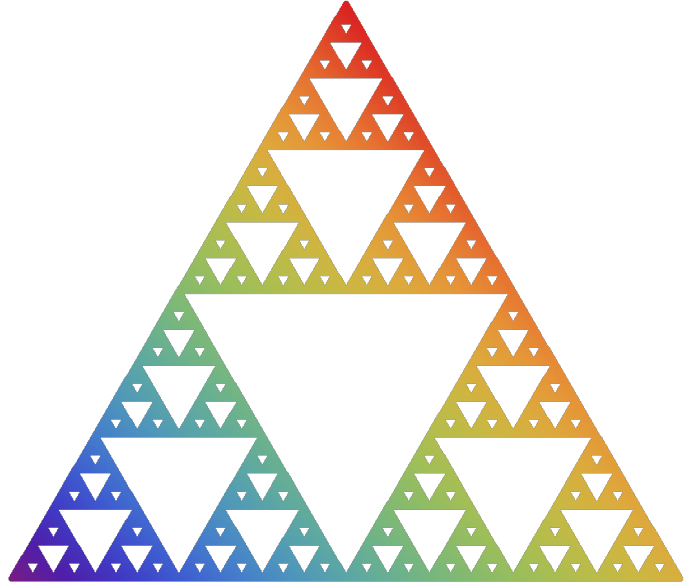}
\captionof{figure}{The graph of the function~$X \mapsto f(X)=\parallel X \parallel^2$.}
\label{fig1}
\end{center}

\begin{center}
\includegraphics[scale=1]{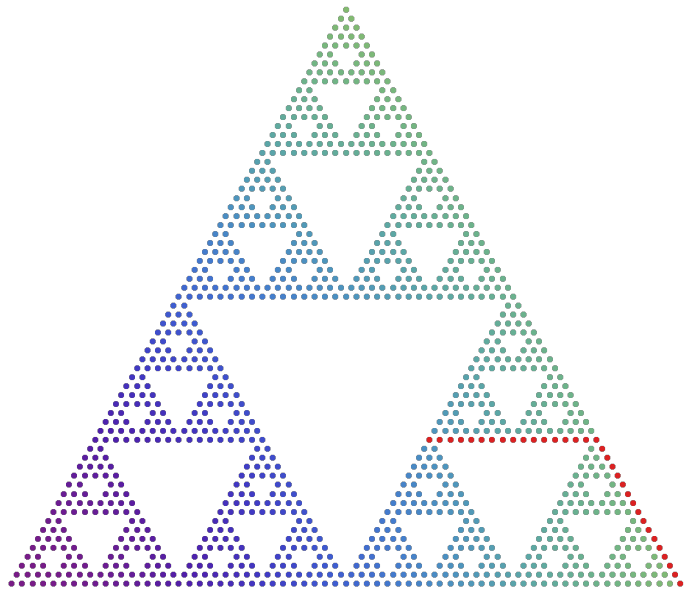}
\captionof{figure}{The algorithm path from $X_0=\displaystyle\left(\frac{5}{8},\frac{\sqrt{3}}{8} \right)$ to $X^{\star}$.}
\label{fig2}
\end{center}

  For this first example, the algorithm starts with $X_0=\displaystyle\left(\frac{5}{8},\frac{\sqrt{3}}{8} \right)$, following the largest gradient (red points), the algorithm converges to the local maximum~$1$ at~$(1,0)$.\\

\begin{center}
\includegraphics[scale=1]{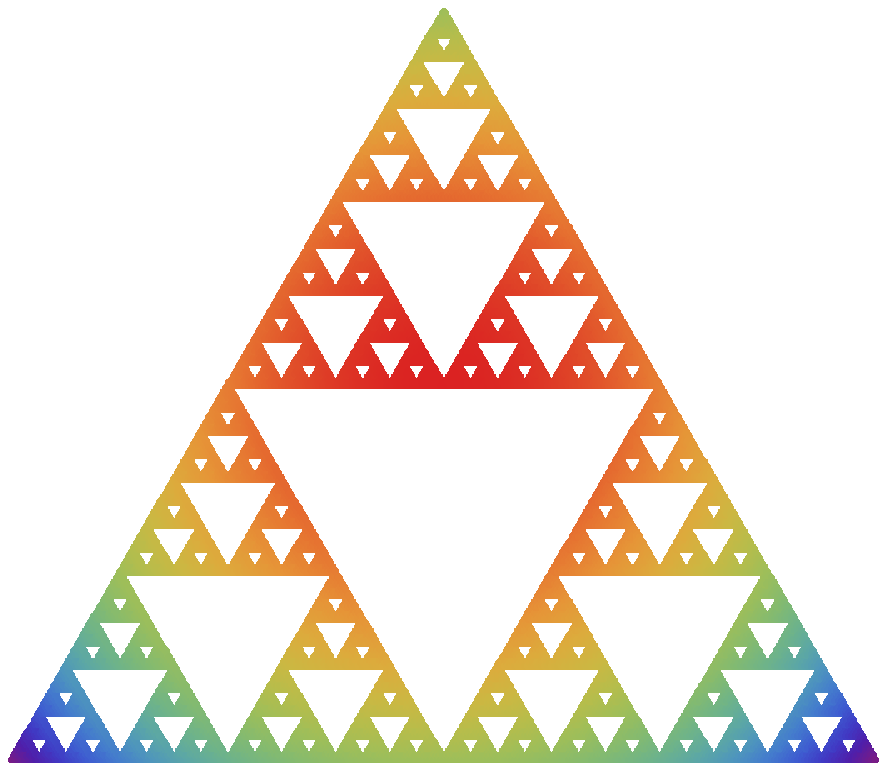}
\captionof{figure}{The graph of the function $X \mapsto g(X)=-\displaystyle\parallel X - \left(\frac{1}{2}, \displaystyle\frac{\sqrt{3}}{4}\right) \parallel^2$.}
\label{fig3}
\end{center}

\begin{center}
\includegraphics[scale=1]{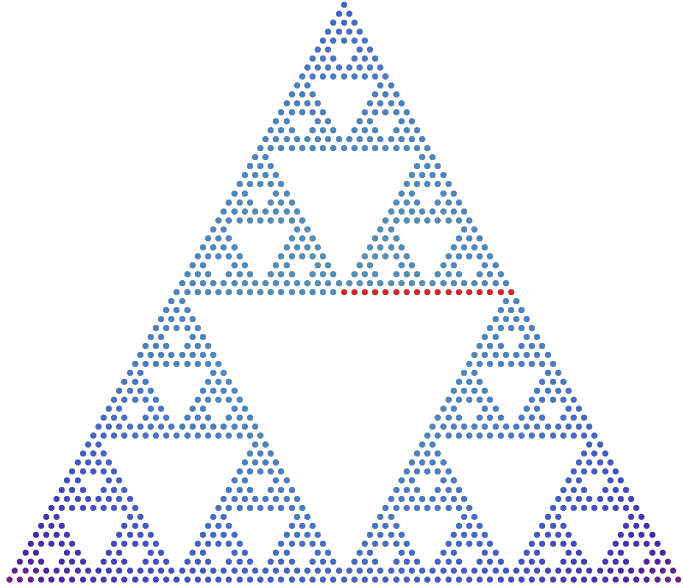}
\captionof{figure}{The algorithm path from $X_0=\displaystyle\left(\frac{3}{4},\frac{\sqrt{3}}{4} \right)$ to $X^{\star}$.}
\label{fig4}
\end{center}

In the second example, the algorithm starts with $X_0=\displaystyle\left(\frac{3}{4},\frac{\sqrt{3}}{4} \right)$, and converges to the global maximum $0$ at $\displaystyle\left(\frac{1}{2},\frac{\sqrt{3}}{4}\right)$.\\
\vskip 1cm

\paragraph{The Sierpi\'{n}ski Tetrahedron \vskip 0.5cm}

The Sierpi\'{n}ski Tetrahedron requires a calculation time for the maximum at step~\mbox{$m\,\in\,\N^\star$} which is of order~$\mathcal{O}(4^{m})$.\\

  In the sequel (see~figures~5 and~6), we present results of our algorithm in the case of Sierpi\'{n}ski Tetrahedron, with vertices:
$$P_0=(0,0,0) \quad , \quad P_1=(1,0,0) \quad , \quad P_2=\left(\frac{1}{2},\frac{\sqrt{3}}{2},0\right)
\quad , \quad P_3=\left(\frac{1}{2},\frac{1}{2\sqrt{3}},\sqrt{\frac{2}{3}}\right)$$

\noindent for the value~$m=6$.\\

 The color function is related to the gradient of the one at stake, high values ranging from red to blue.  

\begin{center}
\includegraphics[scale=1]{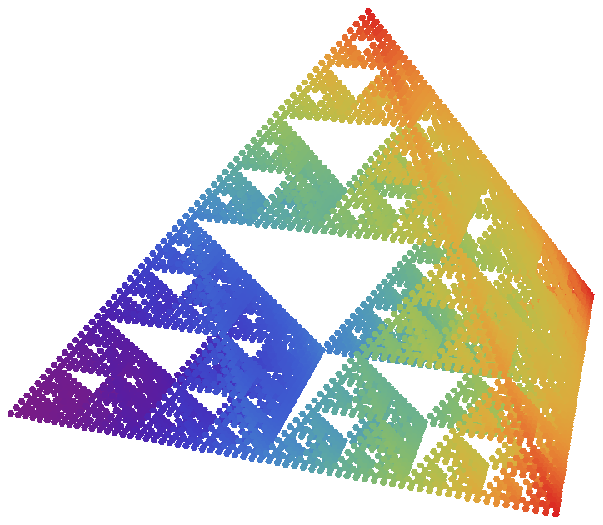}
\captionof{figure}{The graph of the function~$X \mapsto f(X)=\parallel X \parallel^2$.}
\label{fig5}
\end{center}

\begin{center}
\includegraphics[scale=1]{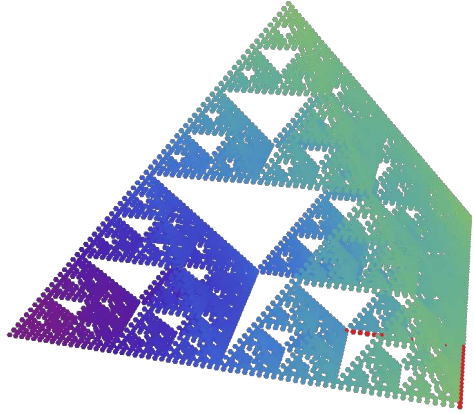}
\captionof{figure}{The algorithm path from $X_0=\displaystyle\left(\frac{41}{64},\frac{1}{64}+\frac{\sqrt{3}}{8},\frac{1}{64}\right)$ to $X^{\star}$.}
\label{fig6}
\end{center}

  For this first example, the algorithm start with $X_0=\displaystyle\left(\frac{41}{64},\frac{1}{64}+\frac{\sqrt{3}}{8},\frac{1}{64}\right)$, following the largest gradient (red points), the algorithm converges to the local maximum~$1$ at $(1,0,0)$.\\
\vskip 1cm

\subsubsection{Self-similar curves}

  Self-similar curves require a calculation time which is of order $\mathcal{O}(N^{m})$, due to the fact that every vertex has only two neighbors. In such cases:~$V_0=\{P_0,P_1\}$. Thus, every vertex~$X$ has exactly two adresses:

$$ X=f_{{\cal W}_i}(P_0)=f_{{\cal W}_{i+1}}(P_1)$$

\noindent where~$\left ({\cal W}_i,{\cal W}_j \right)\, \in \,\Sigma^2$. The neighbors of~$X$ are thus given by:

$$ f_{{\cal W}_{i-1}}(P_0) \quad \text{and}\quad  f_{{\cal W}_{j+1}}(P_1) $$

\noindent where ${\cal W}_{i-1}$ (resp.  ${\cal W}_{j+1}$) is the next (resp. the past) adress of ${\cal W}_i$ (resp. ${\cal W}_j$) in the lexicographical order.\\

  In the sequel (see~figures~7 and~8), we present results of our algorithm in the case of the  Minkowski curve, with $V_0=\left \lbrace (0,0);(1,0) \right \rbrace $, for the value~$m=3$.\\

 The color function is related to the gradient of the one at stake, high values ranging from red to blue.  

\begin{center}
\includegraphics[scale=1]{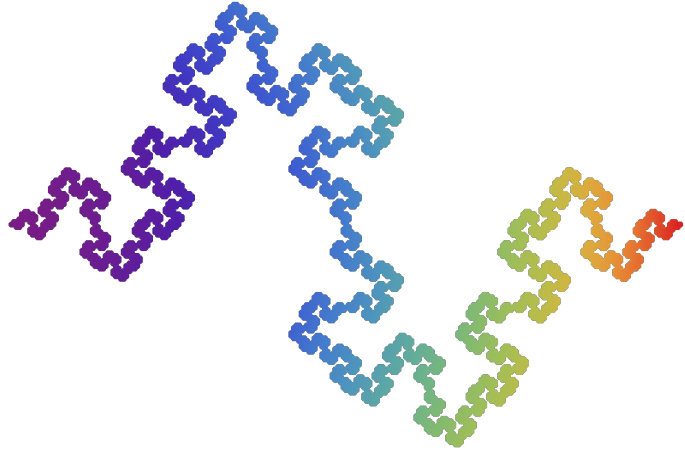}
\captionof{figure}{The graph of the function~$X \mapsto f(X)=\parallel X \parallel^2$.}
\label{fig7}
\end{center}

\begin{center}
\includegraphics[scale=1]{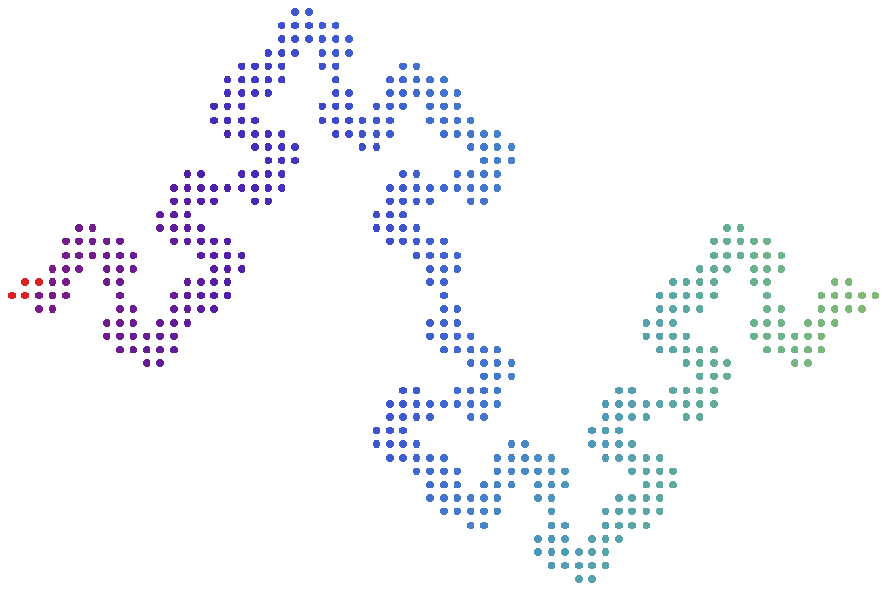}
\captionof{figure}{The algorithm path from $X_0=\displaystyle\left(0,0\right)$ to $X^{\star}$.}
\label{fig8}
\end{center}

\vskip 1cm
  For this example, the algorithm starts with~$X_0=\displaystyle\left(0,0\right)$ ; following the largest gradient (red points), the algorithm converges to the local maximum $\displaystyle\left(\frac{5}{4096}\right)$ at $\displaystyle\left(\frac{1}{32}, \frac{1}{64}\right)$.

\bibliographystyle{alpha}
\bibliography{BibliographieClaire}

\end{document}